\documentclass[letterpaper, 10 pt, conference]{ieeeconf}  

\usepackage{graphicx}
\usepackage{amsmath,amsfonts,mathrsfs,amssymb,dsfont}
\usepackage{booktabs}
\usepackage{mathtools}
\usepackage{ stmaryrd }
\usepackage{cite}
\usepackage{color}
\usepackage{sidecap}
\sidecaptionvpos{figure}{l}
\usepackage{bbm}

\usepackage{bbm}
\newcommand{\et}{\hfill \blacksquare}

\IEEEoverridecommandlockouts                              

\title{\LARGE \bf
Moment dynamics for a class of time-triggered stochastic hybrid systems
}

\author{Mohammad Soltani$^{1}$, Abhyudai Singh$^{2}$
\thanks{$^{1}$M. Soltani is with Department of Electrical and Computer Engineering, University of Delaware, Newark, DE USA 19716.
{\tt\small msoltani@udel.edu}}%
\thanks{$^{2}$A. Singh is with the Department of Electrical and Computer Engineering, Biomedical Engineering, Mathematical Sciences, Center for Bioinformatics and Computational Biology, University of Delaware, Newark, DE USA 19716.
{\tt\small absingh@udel.edu}}}

\begin{document}

\maketitle
\thispagestyle{empty}
\pagestyle{empty}

\begin{abstract}
Stochastic Hybrid Systems (SHS) constitute an important class of mathematical models that integrate discrete stochastic events with continuous dynamics. The time evolution of statistical moments is generally not closed for SHS, in the sense that the time derivative of the lower-order moments depends on higher-order moments. Here, we identify an important class of SHS where  moment dynamics is automatically closed, and hence moments can be computed exactly by solving a system of coupled differential equations. This class is referred to as linear time-triggered SHS (TTSHS), where the state evolves according to a linear dynamical system. Stochastic events occur at discrete times and the intervals between them are independent random variables that follow a general class of probability distributions. Moreover, whenever the event occurs, the state of the SHS changes randomly based on  a probability distribution. Our approach relies on embedding a Markov chain based on phase-type processes to model timing of events, and showing that the resulting system has closed  moment dynamics. Interestingly, we identify a subclass of linear TTSHS, where the first and second-order moments depend only on the mean time interval between events and invariant of their higher-order statistics. Finally, we discuss applicability of our results to different application areas such as network control systems and systems biology.  
\end{abstract}

\section{Introduction}
Stochastic hybrid systems (SHS) are increasingly being used to model noise and uncertainty in physical, biological and engineering systems. Specific applications
include communication networks \cite{HespanhaMar04,BohacekHespanhaLeeObraczkaJun03,jia06}, network control systems \cite{ahs13,hes14},  air traffic control \cite{vgl06,prj09}, biological systems \cite{dsp15,ans14,Singh010,ljs03,bop08,efl06,hls04,rkr09}, power grids \cite{dcd13, wac11}, modeling of energy grids and smart buildings \cite{sma12,adl12}. We refer interested readers to \cite{hes06,tss14,hls00} for a detailed mathematical  characterization of SHS. 

Traditional analysis of SHS relies heavily on various Monte Carlo simulation techniques, which come at a significant computational cost  \cite{hes05,jup09}.  Since one is often interested in computing only the lower-order moments of the state variables, much time and effort can be saved by directly computing these statistical moments without having to run Monte Carlo simulations. Unfortunately, moment calculations in SHS can be non-trivial due to the problem of unclosed dynamics: the time evolution of lower-order moments depends on higher-order moments \cite{sih05}. Closure methods, that approximate higher-order moments as nonlinear functions of lower-order moments, are generally employed in such cases\cite{lkk09,sih10,gil09,svs15,jdd14,sih07ny}. 

The problem of moment closure leads to an interesting question: are there classes of SHS where moments can be computed exactly without the need for closure techniques? Here, we identify such a class of SHS known as 
 time-triggered SHS (TTSHS) that are a special case of piecewise-deterministic Markov processes \cite{dav93,cod08}.  The main ingredients of TTSHS are as follows:
 \begin{enumerate}
\item A continuous state $\bold{x}(t)\in\mathbb{R}^{n}$ that evolves according to a stable linear dynamical system
\begin{equation}
\dot{\bold{x} }(t)= \hat{a}+A  \bold{x}(t),
\label{dynamics0}
\end{equation}
for some constant vector $\hat{a}$ and Hurwitz matrix $A$. 
\item  Timing of events is determined by a renewal process. In particular, stochastic events occur at discrete times $\bold{t}_s, \ s\in \{1,2,\ldots\}$, and the intervals $\bold{t}_s-\bold{t}_{s-1}$ are independent and identical random variables drawn from 
a given probability density function.
\item A reset map defining the change in $\bold{x}$ when the event occurs
\begin{equation}
\bold{x}(\bold{t}_s)\mapsto \bold{x}_+(\bold{t}_s), \label{reset general}
\end{equation}
where $\bold{x}_+(\bold{t}_s)$ denotes the state of the system just after the event. While prior work has considered a deterministic linear reset map
\begin{equation}\label{reset}
\bold{x}(\bold{t}_s)\mapsto J \bold{x}(\bold{t}_s)
\end{equation}
\cite{ahs10,ahs12,ahs13a}, we allow  $\bold{x}_+(\bold{t}_s)$ to be probabilistically determined. \end{enumerate}

Our goal is to connect moments of the continuous state to the statistics of the time interval $\bold{T}\equiv \bold{t}_s-\bold{t}_{s-1}$. The key contribution of this work is to model arrival of events using a phase-type processes  \cite{queueing01}, and show that the resulting systems has closed moment dynamics.
More specifically, the time derivative of an appropriately selected vector of moments depends only on itself, and not on higher-order moments.  As a consequence, moments can be computed exactly by solving a system of differential equations. For the sake of simplicity, we focus on computing the first and second-order moments, but the ideas can be generalized to obtain any higher-order moment.  In addition, a subclass of TTSHS is identified, where the first and second order moments of $\bold{x}$ depend only on the mean time interval between events. In this case, making the arrival of events more random (for fixed mean arrival times) will not result in higher noise in $\bold{x}$. While TTSHS have been considerably used in context of network control systems \cite{hes15}, we illustrate the above methods on an example drawn from cell biology.

\section{TTSHS model formulation}

The timing of events in TTSHS can be modeled through a timer $\boldsymbol\tau$, that measures the time elapsed since the last event. This timer is reset to zero whenever an event occurs and increases over time as $\dot{\boldsymbol\tau}=1$ in between events. Let the time intervals between events follow  a continuous positively-valued
probability density function $f$. Then, the transition intensity for the event is given by the hazard function 
\begin{equation} \label{hr}
h(\tau)\coloneqq  \frac{f(\tau)}{1-\int_{y=0}^{\tau}f(y)dy},
\end{equation}
 \cite{Ross20109,ehp00}. In particular, the probability that an event occurs in the next infinitesimal time interval $(t,t+dt]$  is $h(\tau)dt$. This formulation of the event  arrival process  via 
 a timer allows representation of TTSHS as a state-driven SHS (Fig. 1). Hence, existing tools such as, Kolmogorov equations and Dynkin's formulas for obtaining time evolution of moments can be employed  for studying stochastic dynamics of  TTSHS. Having defined the timing of events, we next focus on how the events alter the state of the system.

\begin{figure}[!thpb]
      \centering
{\includegraphics[width=0.55\columnwidth]{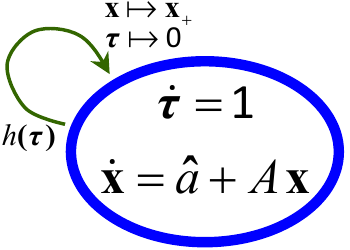}}
      \caption{{\bf Schematic of a linear time-triggered stochastic hybrid system}. The state evolves according to a linear dynamical system and events (green arrow) occur randomly with hazard rate $h(\boldsymbol\tau)$, 
      where the timer $\boldsymbol\tau$ measures the time since the last event. Choosing the hazard rate as \eqref{hr}, ensures that the time between events follows a continuous probability density function $f$. Whenever the event occurs the timer is set to zero and $\bold{x}$ changes via \eqref{reset1}}.
      \label{figmodel}
 \end{figure}

Whenever an event occurs, $\bold{x}$ is reset to $\bold{x}_+$, where $\bold{x}_+$ is a random variable with following statistics  
\begin{subequations}\label{reset1}
\begin{align}
&\langle \bold{x}_+ (\bold{t}_s)\rangle=  J\bold{x}(\bold{t}_s)+ R,\label{conditional x}\\
\left \langle \bold{x}_+(\bold{t}_s)\bold{x}^T_+(\bold{t}_s)\right \rangle& - \langle \bold{x}_+ (\bold{t}_s)\rangle\langle \bold{x}_+ (\bold{t}_s)\rangle^T   = \nonumber
\\& Q\bold{x}(\bold{t}_s)\bold{x}^T(\bold{t}_s) + D \bold{x} (\bold{t}_s)\mathbbm{1}_n+ E,
\label{conditional x2}
\end{align}
\end{subequations}
where the symbol $\langle \ \rangle$ denotes the expected value. Here  $J \in \mathbb{R}^{n\times n}$, $R \in \mathbb{R}^{ n \times 1}$, $ Q \in \mathbb{R}^{n \times n}$, $ D \in \mathbb{R}^{n \times n}$, and $ E \in \mathbb{R}^{ n \times n}$ are constant matrices, $\mathbbm{1}_n$ is a $1 \times n$ unit matrix. Note that the mean of $\bold{x}_+$ is a linear affine function $\bold{x}$, which is a generalization over the linear map \eqref{reset} previously used. The covariation matrix of $\bold{x}_+$ is defined by \eqref{conditional x2} and covers a wide range of possibilities.  For example, $Q=D=E=0$ imply
$\bold{x}_+=J \bold{x}+ R$ with probability one. Moreover, non-zero $Q$, $D$ and $E$ can be used to incorporate constant or state-dependent noise in $\bold{x}_+ $. In the following sections, we 
show how statistical moments of $\bold{x}(t)$ can be computed exactly for TTSHS illustrate in Fig 1. We first consider a subclass of TTSHS, where events only impart noise to the system, in the sense that the average value of $\bold{x}$ just after the event is the same as its value just before the event.

\section{Moment dynamics of TTSHS with noise-imparting events}
Consider a subclass of the TTSHS with $J=I$, $R=Q=0$ reducing \eqref{reset1} to
\begin{subequations}\label{reset2}
\begin{align}
&\langle \bold{x}_+ (\bold{t}_s) \rangle=  \bold{x} (\bold{t}_s),\label{conditional x0}\\
&\left \langle \bold{x}_+\bold{x}_+^T(\bold{t}_s) - \langle \bold{x}_+ (\bold{t}_s) \rangle \langle \bold{x}_+ (\bold{t}_s) \rangle^T \right \rangle = \label{conditional x20}\\ & \hspace{50mm}  D\bold{x} (\bold{t}_s)\mathbbm{1}_n+ E. \nonumber
\end{align}
\end{subequations}
In essence, this reset corresponds to addition of a zero-mean noise term, whenever the event occurs. This noise term can be potentially state-dependent through $D$ in \eqref{conditional x20}.

The time derivative of the expected value of a vector of continuously differentiable functions
\begin{equation}
\varphi(\bold{x},\boldsymbol\tau ) = [\varphi_1(\bold{x},\boldsymbol\tau ), \varphi_2(\bold{x},\boldsymbol\tau ) ,\ldots,  \varphi_n(\bold{x},\boldsymbol\tau )]^T \in \mathbb{R}^n
\end{equation}
can be derived using
\begin{equation}
\begin{aligned}
 \frac{d\langle \varphi(\bold{x},\boldsymbol\tau ) \rangle}{dt}= & \left \langle  ({\textit{L}}\varphi)(\bold{x},\boldsymbol\tau ) \right \rangle, \label{dynnf1}
  \end{aligned}
\end{equation}
where ${\textit{L}}$ is an operator (also called the extended generator of the SHS)
which maps $\varphi$ to $({\textit{L}}\varphi)(x,\tau )$ according to
\begin{equation}
\begin{aligned}
(\textbf{\textit{L}}\varphi)(x,\tau )  \equiv  & \left \langle  \sum_{Events}   h(\tau ) \times \Delta \varphi(x,\tau )  \right \rangle \\
&+ \left \langle \frac{ \partial \varphi(x,\tau )}{\partial x} \left(Ax+\hat{a} \right)   \right \rangle+ \left \langle  \frac{\partial \varphi(x,\tau )}{\partial \tau} \right \rangle, \label{dynnf1}
  \end{aligned}
\end{equation}
\cite{hsi04}. Here, $\Delta \varphi$ is the change in vector $ \varphi$ whenever an event occurs and 
elements of matrix $\frac{ \partial \varphi(x,\tau )}{\partial x}\in \mathbb{R}^{n \times n}$ is defined as
\begin{align}
\left(\frac{ \partial \varphi(x,\tau )}{\partial x}\right)_{ij}= \frac{ \partial \varphi_{i}(x,\tau )}{\partial x_j}.
\end{align}
Setting $ \varphi(\bold{x},\boldsymbol\tau )$ to be $\bold{x}$ in \eqref{dynnf1} and using \eqref{conditional x0}, the time evolution of the 
first-order moments is obtained as 
\begin{equation}
\begin{aligned}
\frac{d\langle \bold{x} \rangle}{dt}
  =& \hat{a} +  A \langle \bold{x} \rangle + \left \langle h(\boldsymbol\tau) \left( \langle \bold{x}_+ \rangle  -\bold{x}\right)\right\rangle = \hat{a} +  A \langle \bold{x} \rangle.
\label{moment dynamics}
\end{aligned}
\end{equation}
For deriving the second-order moments we need the following theorem. \\

\noindent {\bf Theorem 1}:
\emph{ Consider the TTSHS illustrated in Fig. 1 with noise-imparting events, i.e, reset defined by \eqref{reset2}. Then for deterministic initial conditions
\begin{equation}
\begin{aligned}
\langle \bold{x} \vert  \boldsymbol \tau  \rangle=\langle \bold{x}  \rangle
\end{aligned}
\end{equation}
$\forall t > 0$.}$\et$\\
A direct result of this theorem is that 
\begin{equation}
\langle h(\boldsymbol\tau) \bold{x}  \rangle =\langle h(\boldsymbol\tau) \rangle \langle \bold{x}  \rangle .
\end{equation}
Using \eqref{dynnf1} and the above result, the time evolution of the second order moments is obtained as 
\begin{equation}
\begin{aligned}
\frac{d \langle \bold{x}\bold{x}^T \rangle}{dt}  
=&  A  \langle \bold{x}\bold{x}^T \rangle + \langle \bold{x}\bold{x}^T \rangle A^T+\hat{a}\langle \bold{x}^T \rangle + \langle \bold{x} \rangle  \hat{a}^T
\\&+ \langle h(\boldsymbol\tau) (\bold{x}\bold{x}^T+D \bold{x} (\bold{t}_s)\mathbbm{1}_n+ E- \bold{x}\bold{x}^T) \rangle \\
= &  A  \langle \bold{x}\bold{x}^T \rangle +  \langle \bold{x}\bold{x}^T \rangle A^T +\hat{a}\langle \bold{x}^T \rangle  + \langle \bold{x} \rangle  \hat{a}^T \\& + D \langle h(\boldsymbol\tau)\rangle  \langle  \bold{x}  \rangle \mathbbm{1}_n+ E \langle h(\boldsymbol\tau) \rangle.
 \label{second order dunamics}
 \end{aligned}
\end{equation}
At steady-state
\begin{equation}
\begin{aligned}
&A  \overline{\langle \bold{x}\bold{x}^T \rangle }+ \overline{ \langle \bold{x}\bold{x}^T \rangle } A^T +\hat{a}\overline{\langle \bold{x}^T \rangle}  + \overline{ \langle \bold{x} \rangle}  \hat{a}^T= 
\\
&  
- D  \overline{ \langle h(\boldsymbol\tau)  \rangle }\ \overline{\langle  \bold{x}  \rangle }\mathbbm{1}_n
-  \overline{ \langle h(\boldsymbol\tau) \rangle} E ,
 \end{aligned}
 \label{second order dunamics2}
\end{equation}
where $\overline{\langle \ \rangle} $ denotes expected value as $t \to \infty$. Let the random variable $ \bold{T}  \equiv    \bold{t}_s - \bold{t}_{s-1}$ denote the time interval between events for the TTSHS illustrated in Fig. 1.
Then, the mean interval is given by 
\begin{equation}
 \langle  \bold{T} \rangle = \frac{1}{ \overline{ \langle h(\boldsymbol\tau)  \rangle}}  \label{mean cell time}
\end{equation}
 \cite{fin08}.
Using \eqref{mean cell time} and $\overline{\langle  \bold{x}  \rangle}=-A^{-1}\hat{a}$ from \eqref{moment dynamics}, the above equation reduces to
\begin{equation}\label{00}
\begin{aligned}
&A C + C A^T = 
-\frac{1}{ \langle  \bold{T}  \rangle }  D   \overline{\langle  \bold{x}  \rangle }\mathbbm{1}_n
- \frac{1}{ \langle  \bold{T}  \rangle } E,
 \end{aligned}
\end{equation}
where the steady-state covariance matrix
\begin{equation}
\begin{aligned}
C= \overline{\langle \bold{x}\bold{x}^T \rangle} -\overline{\langle \bold{x}\rangle\langle\bold{x}^T \rangle}
 \end{aligned}
\end{equation}
is a unique solution to the Lyapunov equation \eqref{00}. Interestingly, $C$ only depend on $ \langle  \bold{T} \rangle$ and independent of higher-order statistics of $\bold{T} $. Thus, making timing of events more random for a fixed $ \langle  \bold{T} \rangle$ will not affect noise levels in $\bold{x}$. It is important to point out that this result only holds for steady-state, and the transient covariance matrix will depend on the entire distribution of $\bold{T} $.
Next, we illustrate these result on an example drawn from systems biology.

\begin{figure*}[th]
\centering
\vspace{8mm}
\includegraphics[width=0.8\textwidth]{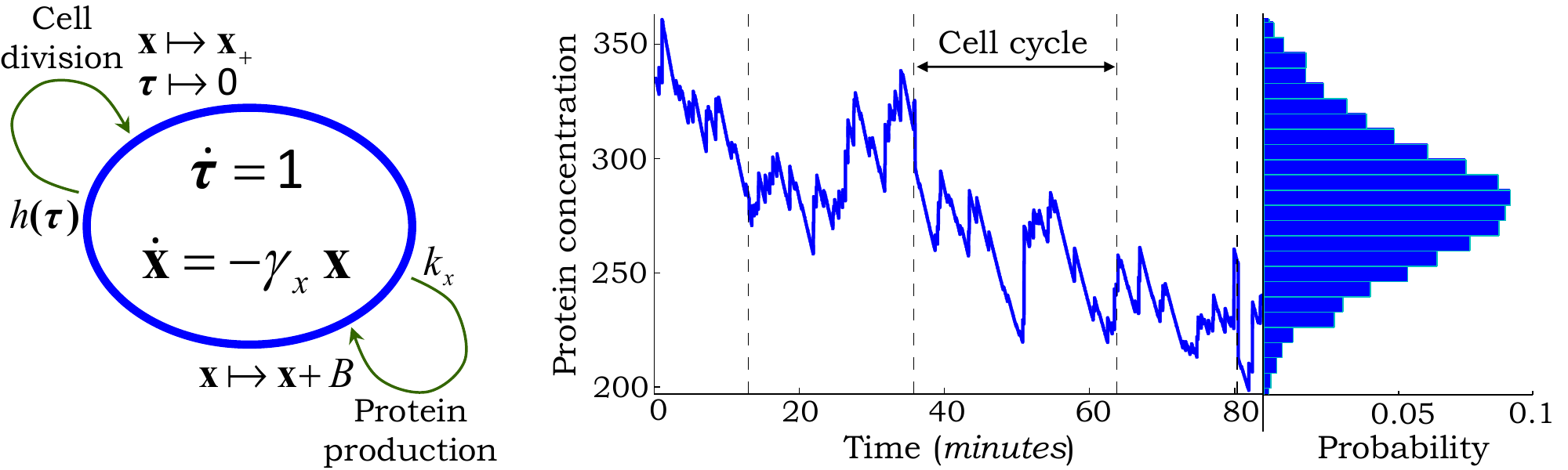}
\caption{
{\bf Modeling stochasticity in the level of protein using TTSHS}. \textit{Left}: Time evolution of the protein level $\bold{x} \in \mathbb{R}$ in a single cell is modeled via a SHS with two stochastic resets representing production of proteins in bursts and cell-division events. The latter is controlled by a timer $\boldsymbol\tau$, and whenever it occurs the state is reset via \eqref{division character000}-\eqref{division character111}. In between events, protein concentrations decay exponentially with rate $\gamma_x$ due to cellular growth. \textit{Right}: A sample trajectory of $\bold{x}$ is shown with cell division events (dashed lines). The steady-state distribution of $\bold{x}$ obtained via a large number Monte Carlo simulations is shown on the right.}
\label{fig2}
\end{figure*}

\section{Systems biology example}

In the previous section, we presented results on moment dynamics of TTSHS with noise-imparting events. Here we consider an example of such a system drawn from cellular biology. 
Advances in experimental technologies have shown  that synthesis of protein molecules from its corresponding gene is an inherently random process \cite{bpm06,rao05}. Such stochasticity results in random fluctuations in the level of a protein inside an individual cell, and  critically impacts functioning of intracellular biological pathways \cite{lod08,vsk08,wds08,siw09,sid13,bos15, sbf15}.  Below we present and analyze a TTSHS model that captures different sources of noise in protein levels.

Let $\bold{x}(t) \in \mathbb{R}$ denote the concentration of a given protein inside a single cell at time $t$. We model the time evolution of $\bold{x}$ using a TTSHS that incorporates three noise mechanisms:
 \begin{enumerate}
\item Stochastic production of protein molecules in bursts of gene activity, as has been experimentally seen \cite{rpt06,src10,smg11}.

\item Timer controlled cell-division events occur at random times. Whenever, cell division occurs,  both the cell volume and the number of protein molecules reduce by half (assuming symmetric division). Thus, in the sense of average concentrations, there is no change
\begin{equation}
 \langle \bold{x}_+(\bold{t}_s) \rangle= \bold{x} (\bold{t}_s), \label{division character000}
\end{equation}

\item During cell division, protein molecules are partitioned between two daughters cells based on a binomial distribution, i.e., each molecule has an equal probability of being in one of the two cells \cite{gpz05,sva15}. This binomial partitioning process introduces noise in the concentration that can be represented by
\begin{equation}
 \left \langle \bold{x}_+^2(\bold{t}_s) - \langle \bold{x}_+(\bold{t}_s) \rangle^2  \right \rangle =   \beta \bold{x}   (\bold{t}_s) .
  \label{division character111}
\end{equation}
for some positive constant $\beta$. The linear dependence in \eqref{division character111} comes from the fact that the variance of a binomially distributed random variable is proportion to its mean.
\end{enumerate}

Before considering the stochastic case, we first consider deterministic protein production, where the concentration evolves as  
\begin{equation}
\dot{\bold{x}}(t)=k_x \langle B \rangle - \gamma_x \bold{x}(t). \label{concentration}
\end{equation}
Here $k_x$ and $\langle B \rangle$ denote the frequency and size of protein bursts, resulting in a net production rate of $k_x \langle B \rangle$ in the deterministic model. The concentration is diluted at a rate $\gamma_x$, which is the rate of exponential growth in cell volume. Cell-division events are assumed to occur randomly at times $\bold{t}_s$, where
$\bold{T}\equiv \bold{t}_s-\bold{t}_{s-1}$ follows an arbitrary positively-valued distribution. The mean cell-division time is intimately  connected to the growth rate via
\begin{equation}
\langle \bold{T} \rangle = \frac{ln(2)}{2\gamma_x}\label{gamma}
\end{equation}
\cite{alo11}. The TTSHS defined by  \eqref{division character000}-\eqref{concentration} is of the form discussed in the Section III, and time evolution of moments is obtained as
\begin{subequations}
\begin{align}
&\frac{d \langle \bold{x} \rangle }{dt} = k_x \langle B \rangle -\gamma_x \langle \bold{x} \rangle, \\
&\frac{d \langle \bold{x}^2 \rangle}{dt}=  2 k_x \langle B \rangle\langle \bold{x} \rangle - 2 \gamma_x \langle \bold{x}^2 \rangle + \beta  \langle h(\boldsymbol\tau) \rangle \langle \bold{x} \rangle.
\end{align}
\end{subequations}
Steady-state analysis yields the following mean and noise (as quantified by the coefficient of variation squared)
\begin{align}
&\overline{\langle \bold{x} \rangle}=\frac{k_x \langle B \rangle  }{\gamma_x}, \ \ \ \ 
&CV^2_x = \frac{\overline{\langle \bold{x}^2 \rangle} - \overline{\langle \bold{x} \rangle} ^2}{\overline{\langle \bold{x} \rangle}^2}=\frac{ln(2)\beta}{2\overline{\langle \bold{x} \rangle}},\label{fano factor concent}
\end{align}
respectively. As predicted by theory, the noise in the protein level only depends on $\langle \bold{T} \rangle$, which enters the equation via $\gamma_x$ (see \eqref{gamma}). Thus remarkably, making the timing of cell division more random (for a fixed mean) will not result in higher stochasticity in the protein concentration. 

Next, we consider stochastic production of proteins, which involves adding a second family of resets in the above TTSHS model. As pointed earlier, production of proteins is assumed to occur in random bursts that happen at random times. In particular, burst events are assumed to occur at exponentially distributed times with rate $k_x$. Whenever the event occurs the concentration changes as
\begin{equation}
\bold{x}(\bold{t}_s) \mapsto \bold{x}(\bold{t}_s)+ B,
\end{equation}
where $B$ is a random variable denoting the size of protein bursts and follows an arbitrary positively-valued distribution. The overall model with both stochastic production and cell division events is shown in Fig. 2. In between events, the concentration is diluted as  
\begin{equation}
\dot{\bold{x}}(t)= - \gamma_x \bold{x}(t). 
\end{equation}
The time evolution of moments for this new SHS is given by the corresponding extended generator and given by  
\begin{subequations}
\begin{align}
&\frac{d \langle \bold{x} \rangle }{dt} = k_x \langle B \rangle -\gamma_x \langle \bold{x} \rangle, \\
&\frac{d \langle \bold{x}^2 \rangle}{dt}=k_x \langle B^2 \rangle+  2 k_x \langle B \rangle\langle \bold{x} \rangle - 2 \gamma_x \langle \bold{x}^2 \rangle + \beta  \langle h(\boldsymbol\tau) \rangle \langle \bold{x} \rangle
\end{align}
\end{subequations}
Steady-state analysis of the above moment equations yields 
\begin{align}
&\overline{\langle \bold{x} \rangle}=\frac{k_x \langle B \rangle  }{\gamma_x}, \ \ \ \ 
&CV^2_x =\frac{ln(2)\beta}{2\overline{\langle \bold{x} \rangle}} +\frac{1}{2}\frac{\langle B^2 \rangle}{\langle B \rangle\overline{\langle \bold{x} \rangle} }, \label{fano factor concent}
\end{align}
providing the first results  connecting the protein noise level to randomness in the underlying bursty synthesis and cell division events. 
As before, the noise only depends on $\langle \bold{T} \rangle$ and independent of its higher-order statistics. 

\begin{table*}[!t]
\centering
\vspace{5mm}
\caption{Stochastic events in the TTSHS with timing of events described by the embedded Markov Chain in Fig. 3.  }
\label{table reaction}
\begin{center}
\begin{tabular}{ccc}
\toprule
Stochastic events & Reset & Transitions intensity \\
 \midrule
Phase-type evolution  &$\bold{s}_{ij}(t)\mapsto \bold{s}_{ij}(t) -1 $, $ \ \ \bold{s}_{i(j+1)}(t)\mapsto \bold{s}_{i(j+1)}(t) +1$&$k_i \bold{s}_{ij}$, $\ \ j \in \{ 1,\ldots,m_i-1 \}$ \\ 
\\
Events changing $\bold{x}$ &$\bold{x}(\bold{t}_s)\mapsto \bold{x}_+(\bold{t}_s), $ $\ \ \bold{s}_{jm_j}(\bold{t}_s)\mapsto 0, $ $ \ \ \bold{s}_{i1}(\bold{t}_s)\mapsto \bold{s}_{i1}(\bold{t}_s) +1 $&$p_{i} \sum_{j=1}^{n} k_j \bold{s}_{j{m_j}}$, $ \ \ i \& j \in \{1,\ldots, n \}$\\
\bottomrule
\end{tabular}
\end{center}
\end{table*}

\section{Moment dynamics of TTSHS}

In the previous sections, we have considered a sub-class of TTSHS with resets given by \eqref{reset2} (i.e., noise-imparting events). Here, we return to the original TTSHS described by 
\eqref{dynamics0} \& \eqref{reset1}, and show how time evolution of moments can be computed for these systems.  

In general, moment equations cannot be solved exactly for any arbitrary function $h(\boldsymbol\tau)$ in Fig. 1. Our strategy for exact moment computations relies on two steps: i) Modeling the timing of stochastic events through a phase-type distribution, which can be represented by embedding a continuous-time Markov chain (Fig. 3) and ii)
Showing that the time evolution of moments in the resulting system becomes automatically closed at some high-order moment. Here we focus on phase-type distributions that consists of a mixture of Erlang distributions \cite{queueing01}, and use them as a practical tool for modeling the timing of stochastic events in TTSHS.

{\bf Embedded Markov chain for event timing}: Recall that an Erlang distribution of order $i$ is the distribution of the sum of $i$ independent and identical exponential random variables.  The interval $\textbf{T}\equiv \bold{t}_s-\bold{t}_{s-1}$ is assumed to have an Erlang distribution of order $m_i$ with probability $p_i, \ i=\{1,\ldots,n\}$ and can be represented by a continuous-time Markov chain with states $S_{ij}$, $j=\{1,\ldots,m_i\}, \ i=\{1,\ldots,n\}$ (Fig. 3). Let Bernoulli random variables $\bold{s}_{ij}=1$ if the system resides in state $S_{ij}$ and $0$ otherwise.  The probability of transition $S_{ij}\rightarrow S_{i(j+1)}$ in the next infinitesimal time interval 
$(t,t+dt]$ is given by $k_i\bold{s}_{ij}dt$, implying that the time spent in each state $S_{ij}$ is  exponentially distributed with mean $1/k_i$. To summarize, just after an event occurs a state $S_{i1}, \ i=\{1,\ldots,n\}$ is chosen with probability $p_i$ and the next event occurs after transitioning through $m_i$ exponentially distributed steps. Based on this formulation, the probability of the stochastic event occurring in the time interval $(t,t+dt]$ is given by $  p_i \sum_{j=1}^n k_j\bold{s}_{jj}dt$, and whenever the event occurs, the state is reset as per \eqref{reset1}.  For a mixture of Erlang distributions, the moment are  given by 
\begin{equation}
\begin{aligned}
& \langle \bold{T}^q \rangle = \sum_{i=1}^{n}\frac{p_i}{k_i^q}\frac{(m_i+q-1)!}{(m_i-1)!} , \ \ q \in \{1,2,\ldots\}
\label{phase mean cv}
\end{aligned}
\end{equation}
\cite{bkf14}. Given a specific distribution of timing of events, the above equation can in principle be used to construct a complex enough appropriate Marko chain that matches some lower orders moments of the given distribution.

\begin{figure}[!h]
\includegraphics[width=0.5\textwidth]{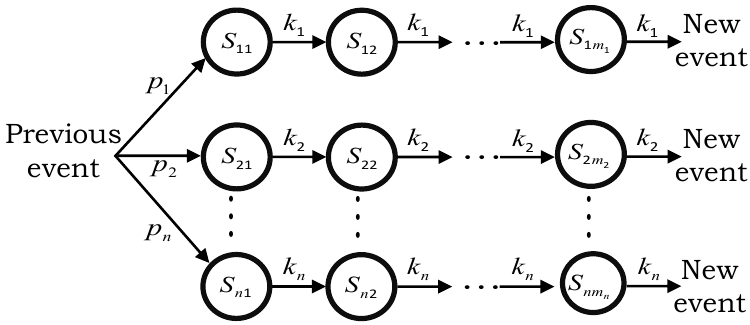} 
\caption{
{\bf A continuous-time Markov chain model for timing of events in TTSHS}. The time interval $T\equiv \bold{t}_s-\bold{t}_{s-1}$ between two successive stochastic events is assumed to follow a mixture of Erlang distributions. 
After an event occurs, a state $S_{i1}$, $i=\{1,\ldots,n\}$ is chosen with probability $p_i$. The systems transitions through states $S_{ij}, \ j=\{1,\ldots,m_i\}$ residing for an exponentially distributed time in each state. The next event occurs after exit from $S_{im_i}$ and the above process is repeated.
}
\label{fig3}
\end{figure}
The overall model is now given by the linear system
\begin{equation}
\dot{\bold{x}}(t) = \hat{a} + A \bold{x}. \label{dynamics02}
\end{equation}
together with stochastic transitions associated with the embedded Markov chain illustrated in  Table \ref{table reaction}. The theorem below outlines the main result.\\

\noindent {\bf Theorem 2}:
\emph{Consider the TTSHS where timing of events are modeled through the Markov chain in Fig. 3, and reset are given by \eqref{reset1}. Then, there exists a vector 
 $\mu$ of all first and second order moments of stochastic processes $\bold{x}$ and $\bold{s}_{ij}$, and selected third order moments, such that its time evolution is given by a linear dynamical system
\begin{equation}
\label{moments}
\small
 {{\dot{{\mu}}}}={a_1}+A_1 {{{\mu}}},
\end{equation}
for an appropriate vector $a_1$ and matrix $A_1$.} $\et$\\
\\{\bf Proof of Theorem 2}:
Time derivative of the expected value of any vector of continuously differentiable functions
$\varphi(\bold{x},\bold{s}_{ij})$ is given by 
\begin{equation}
\begin{aligned}
\frac{d\langle \varphi(\bold{x},\bold{s}_{ij}) \rangle}{dt}= \left \langle  ({\textit{L}}\varphi)(\bold{x},\bold{s}_{ij}) \right \rangle,   
\label{dynnf}
\end{aligned}
\end{equation}
where the extended generator $({\textit{L}}\varphi)$ for this SHS is 
\begin{equation}
\begin{aligned}
 ({\textit{L}}\varphi)(x,s_{ij}) =& \left \langle \sum_{Events}  \Delta \varphi(x,s_{ij}) \times \psi (x,s_{ij})\right \rangle \\
&+ \left \langle \frac{\partial \varphi(x,s_{ij})}{\partial x} \left(\hat{a}+Ax \right)\right \rangle. 
\label{dynnf}
\end{aligned}
\end{equation}
Here  $ \psi (x,s_{ij})$ denotes the transition intensities for the events and determine how often these events occur \cite{hsi04}. 

Using the   transition intensities shown in Table I, the dynamics of the means can be written as
\begin{small}
\begin{subequations}
\begin{align}
&\frac{d\langle \bold{x} \rangle}{dt}= \hat{a} + A\bold{x}  + \sum_{j=1}^n  k_j \left( (J-I)\langle\bold{x} \bold{s}_{jm_j} \rangle + R\langle \bold{s}_{jm_j} \rangle \right),\\
&\frac{d\langle \bold{s}_{i1} \rangle}{dt} =  p_i\left \langle \sum_{j=1}^n k_j \bold{s}_{jm_j} \right\rangle - k_i \langle \bold{s}_{i1} \rangle ,\  \ \ i= \lbrace 1,\ldots,n \rbrace, \\
&\frac{d\langle \bold{s}_{ij} \rangle}{dt} =k_i \langle \bold{s}_{i(j-1)} \rangle - k_i \langle \bold{s}_{ij} \rangle ,\\
& \hspace{29mm} i= \lbrace 1,\ldots,n \rbrace,\ \
 j= \lbrace 2,\ldots,m_i \rbrace. \nonumber
\end{align}
\label{mean01} 
\end{subequations}
\end{small}
Note the the first equation is not closed since it depends on the second order moments of the form $\langle\bold{x} \bold{s}_{ij} \rangle$. The time evolution of the moments 
$\langle\bold{x} \bold{s}_{ij} \rangle$ depends on third order moments of the form $\langle  \bold{x} \bold{s}_{ij}^2   \rangle$ and $\langle \bold{s}_{ij} \bold{s}_{rq} \bold{x}^b \rangle$.
However, using the fact that $\bold{s}_{ij}$ are Bernoulli random variables
\begin{equation}\label{22}
\begin{aligned}
\langle \bold{s}_{ij}^q \rangle=\langle \bold{s}_{ij}\rangle,  \ \ 
\langle \bold{s}_{ij}^q \bold{x}^b \rangle=\langle \bold{s}_{ij} \bold{x}^b \rangle, \ \ q \ \& \ b\in \{1,2,\ldots\}.
\end{aligned}
\end{equation} 
Moreover, since only one of the states $\bold{s}_{ij}$ can be $1$ at a time
\begin{equation}\label{24}
\langle \bold{s}_{ij} \bold{s}_{rq} \bold{x}^b \rangle=0, \ {\rm if} \ i\neq r \ {\rm or} \ j\neq q.
\end{equation}
Exploiting \eqref{22}-\eqref{24}, moment dynamics of $\langle \bold{x} \bold{s}_{ij} \rangle$ becomes automatically closed and given by 
\begin{subequations}
\begin{align}
& \frac{d\langle\bold{x} \bold{s}_{i1} \rangle}{dt}= \hat{a}   \langle\bold{x} \bold{s}_{i1} \rangle + A \langle\bold{x} \bold{s}_{i1}\rangle-  k_i\langle \bold{x}\bold{s}_{i1} \rangle \label{final xs}
\\& +  p_i\sum_{j=1}^n k_j \left( (J-I) \left \langle  \bold{x}\bold{s}_{jm_j}\right\rangle  + R  \langle  \bold{s}_{jm_j} \rangle \right)  , \
 i= \lbrace 1,\ldots,n \rbrace, \nonumber \\
& \frac{d\langle \bold{x} \bold{s}_{ij} \rangle}{dt}=  \hat{a}   \langle\bold{x} \bold{s}_{ij} \rangle + A \langle\bold{x} \bold{s}_{ij}\rangle  
-  k_i\langle \bold{x} \bold{s}_{ij} \rangle  \label{gene total2}  \\
& \hspace{11mm}+  k_i \langle \bold{x}  \bold{s}_{i(j-1)} \rangle, \  \
 i= \lbrace 1,\ldots,n \rbrace, \ \
 j= \lbrace 2,\ldots,m_i \rbrace.\nonumber
\end{align}
\label{xsij} 
\end{subequations}
Thus \eqref{mean01} and \eqref{xsij} represent a closed set of equations that can be used to obtain the mean dynamics $\langle \bold{x} \rangle$. A similar approach can be taken for obtaining the second order moments. Briefly, the time evolution of  $\langle \bold{xx}^T \rangle$ would depend on moments of the form $\langle \bold{xx}^T \bold{s}_{ij}\rangle$. Moment dynamics of $\langle \bold{xx}^T \bold{s}_{ij}\rangle$ can be closed automatically using \eqref{22}-\eqref{24}.
$\et$ \\

Given space limitations, we do not present any example illustrating the results of Theorem 2. In summary, our results show that if timing events in TTSHS of can be modeled via a phase-type process (as in Fig. 3), then time evolution of moments can be obtained by solving a linear systems of differential equations.

\section{Conclusion}

Here we studied stochastic dynamics for a class of piecewise-deterministic Markov processes, known as  TTSHS (Fig. 1). In essence, TTSHS are linear dynamical systems, whose state is initialized by random initial conditions at random times. The main contribution of this study is to show how moments can be obtained exactly by solving a closed system of differential equations. Our study also identifies a subclass of TTSHS given by \eqref{reset2}, where moments of $\bold{x}$ depend only on the mean time interval between events. This result leads to an intriguing finding in the context of cell biology: the magnitude of noise in the 
concentration of a protein is invariant of the noise in cell-division times. We are currently working closely with experimental collaborators to test this prediction in living cells. 

It is important to point out that the results presented here can be easily generalized to consider the following scenarios:
\begin{enumerate}
\item The state evolves according to a linear stochastic differential equation between events. 
\item Multiple family of resets, where some resets are timer-dependent, while other are state-dependent. In the latter case, the probability of event occurrence is a linear affine function of $\bold{x}$.
\end{enumerate}
The second point is exemplified from the example shown in Fig. 2 with two families of stochastic resets, albeit the timing of one of the resets was memoryless.

Our study also presents exciting avenues for future research. One direction of research would be to find time evolution of moments for TTSHS with multiple discrete modes, allowing for switching between linear dynamical systems (Fig. 4). The current formulation of TTSHS considers time intervals between events to be independent. It will be interesting to add some form of correlation between successive events. This is particularly important for cell division, where the cell-cycle lengths of mother and daughter cells are generally correlated. Finally, recent work has identified classes of nonlinear stochastic systems, where moment dynamics becomes automatically closed at some higher-order moment \cite{sos15}. In light of this study, one can explore  nonlinearities in TTSHS such that moments can be computed without requiring closure schemes.

\begin{figure}[!h]
\centering
\vspace{5mm}
\includegraphics[width=0.46\textwidth]{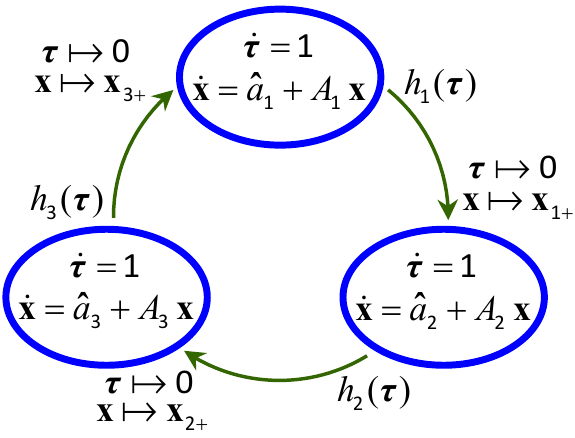} 
\caption{
{\bf Schematic of a TTSHS with multiple discrete modes}. In each discrete model, the state evolves according a linear dynamical system. Transition between discrete states occur based on hazard functions $h_i(\boldsymbol\tau), i=\{1,2,\ldots \}$, and whenever the transition occurs, the state is reset based on maps similar to \eqref{reset1}.
}
\label{fig4}
\end{figure}

\section*{ACKNOWLEDGMENT}
AS is supported by the National Science Foundation Grant DMS-1312926.

\bibliographystyle{IEEEtran}       
\bibliography{references}

\end{document}